\documentclass[preprint]{elsarticle}

\usepackage{amssymb,amsmath,amsthm}
\usepackage{graphics, graphicx, epsfig}

\newtheorem{theorem}{Theorem}[section]

\makeatletter
\def\ps@pprintTitle{%
  \let\@oddhead\@empty
  \let\@evenhead\@empty
  \let\@oddfoot\@empty
  \let\@evenfoot\@oddfoot
}
\makeatother

\begin{document}

\title{Spanning trees with small degrees and few leaves\tnoteref{label1}}
\tnotetext[label1]{Research supported by Conacyt, M\'exico}

\date {}

\author {Eduardo Rivera-Campo}
\ead{erc@xanum.uam.mx}

\address {Departamento de Matem\'aticas \\ Universidad Aut\'onoma Metropolitana -
Iztapalapa \\ Av. San Rafael Atlixco 186, M\'exico D.F. 09340,
M\'exico}


\begin{abstract}

We give an Ore-type condition sufficient for a graph $G$ to have a
spanning tree with small degrees and with few leaves.
\end{abstract}


\maketitle


\section{Introduction}

From a classical result by Ore \cite{O} it is well-known that if a
simple graph $G$ with $n\geq 2$ vertices is such that $d\left(
u\right) +d\left( v\right) \geq n-1$ for each pair $u,v$ of
non-adjacent vertices of $G$, then $G$ contains a hamiltonian path.

A \emph{leaf }of a tree $T$ is a vertex of $T$ with degree one. A
natural generalisation of hamiltonian paths are spanning trees with
a small number of leaves. In this direction, Ore's result was
generalised by Broersma and Tuinstra \cite{BT} \ to the following
theorem.

\begin{theorem}
$\emph{[1]}$ Let $s\geq 2$ and $n\geq 2$ be integers. If $G$ is a
connected simple graph with $n$ vertices such that $d\left( u\right)
+d\left( v\right) \geq n-s+1$, for each pair $u,v$ of non-adjacent
vertices, then $G$ contains a spanning tree with at most $s$ leaves.
\end{theorem}

Further related results have been obtained by Egawa \textit{et al}
\cite {EMYY} and by Tsugaki and Yamashita \cite{TY}. See also
\cite{OY} for a survey on spanning trees with specific properties.

In this note we consider spanning trees with small degrees as well
as with a small number of leaves. Our result is the following.

\begin{theorem}
Let $n$, $k$ and $d_{1},d_{2},\ldots ,d_{n}$ be integers with $1\leq
k\leq n-1$ and $2\leq d_{1}\leq d_{2}\leq \cdots \leq d_{n}\leq
n-1$. If $G$ is a $ k$-connected simple graph with vertex set
$V\left( G\right) =\left\{ w_{1},w_{2},\ldots ,w_{n}\right\} $ such
that $d\left( u\right) +d\left( v\right) \geq
n-1-\sum\limits_{j=1}^{k}\left( d_{i}-2\right) $ for any
non-adjacent vertices $u$ and $v$ of $G$, then $G$ has a spanning
tree $T$ with at most $2+\sum\limits_{j=1}^{k}\left( d_{j}-2\right)
$ leaves and such that $d_{T}\left( w_{j}\right) \leq d_{j}$ for
$j=1,2,\ldots ,n$.
\end{theorem}


\section{Proof of Theorem 2}

Let $T$ be a largest subtree of $G$ with at most
$2+\sum\limits_{j=1}^{k} \left( d_{j}-2\right) $ leaves and such
that if $w_{j}\in V\left( T\right) $ , then $d_{T}\left(
w_{j}\right) \leq d_{j}$. Since $G$ is $k$-connected and $n\geq 2$,
it contains a path with at least $k+1$ vertices. Therefore, we may
assume that tree $T$ has at least $k+1$ vertices.

If $T$ is not a spanning tree, there is a vertex $w$ of $G$ not in
$T$. By Menger's theorem, there are $k$ internally disjoint paths
$\pi _{1},\pi _{2},\ldots ,\pi _{k}$ in $G$ joining $w$ to $k$
different vertices $r_{1},r_{2},\ldots ,r_{k}$ of $T$.

Let $n_{1}$ denote the number of leaves of $T$. We claim
$n_{1}=2+\sum\limits_{j=1}^{k}\left( d_{j}-2\right) $, otherwise
there is a vertex $r_{i}$ such that $d_{T}\left( r_{i}\right)
<d_{j_{i}}$ where $ w_{j_{i}}=r_{i}$. Then $T^{\prime }=T\cup \pi
_{i}$ is a subtree of $G$ with more vertices than $T$ such that
$d_{T^{\prime }}\left( w_{j}\right) \leq d_{j}$ for each $w_{j}\in
V\left( T^{\prime }\right) $ and with at most $ n_{1}+1\leq
2+\sum\limits_{j=1}^{k}\left( d_{j}-2\right) $ leaves, which
contradicts our assumption on the maximality of $T$.

Because of Ore's theorem, we can assume $d_{i}\geq 3$ for some
$i=1,2,\ldots ,k$. Since $T$ has
$n_{1}=2+\sum\limits_{j=1}^{k}\left( d_{j}-2\right) \geq 3$ leaves,
as shown above, there is a vertex $w_{j}$ of $T$ such that $
d_{T}\left( w_{j}\right) \geq 3$. Suppose there are vertices $x$ and
$y$ of degree one in $T$ such that $xy\in E\left( G\right) $. Since
$T$ is not a path, there is an edge $zz^{\prime }$ in the unique
$xy$ path contained in $ T $ with $d_{T}\left( z\right) \geq 3$. Let
$T^{\prime }=\left( T-zz^{\prime }\right) +xy$ and notice that
$T^{\prime }$ is a subtree of $G$ with $
V\left( T^{\prime }\right) =V\left( T\right) $, with less than $
2+\sum\limits_{j=1}^{k}\left( d_{j}-2\right) $ leaves and such that
$ d_{T^{\prime }}\left( w_{j}\right) \leq d_{j}$ for each $w_{j}\in
V\left( T^{\prime }\right) $. As above, this is a contradiction and
therefore no leaves of $T$ are adjacent in $G$.

\begin{figure}[h!!]
\begin{center}
\includegraphics[height=4cm, width=5cm]{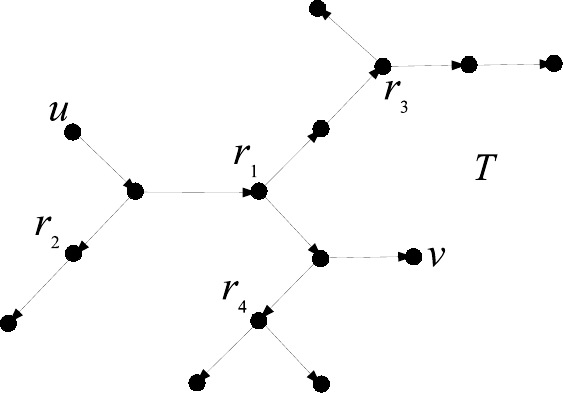}
\caption{$n=15$, $k=4$, $d_1,d_2,\ldots, d_{15}=3$.} \label{one}
\end{center}
\end{figure}

Notice that $d_{T}\left( r_{1}\right) \geq 2$, otherwise $T^{\prime
}=T\cup \pi _{1}$ would be a tree larger than $T$, with the same
number of leaves and with $d_{T^{\prime }}\left( w_{j}\right) \leq
d_{j}$ for each vertex $ w_{j}$ of $T^{\prime }$. Let $u$ and $v$ be
any two leaves of $T$ with the property that the vertex $r_{1}$ lies
in the unique $uv$ path $T_{uv}$, contained in $T$. Orient the edges
of $T$ in such a way that the corresponding directed tree
$\overrightarrow{T}$ is outdirected with root $u$ (see Fig.
\ref{one}.)

For each vertex $z\neq u$ in $T$ let $z^{-}$ be the unique vertex of
$T$ such that $z^{-}z$ is an arc of $\overrightarrow{T}$. Let
\[
A=\left\{ y\in V\left( T\right) :yv\in E\left( G\right) \right\}
\text{ and } B=\left\{ x^{-}\in V\left( T\right) :ux\in E\left(
G\right) \right\} .
\]

Because of the way the tree $T$ was chosen, all vertices of $G$
adjacent to $ u$ or to $v$ lie in $T$ and therefore $\left\vert
A\right\vert =d\left( v\right) $. Let $x_{1}$ and $x_{2}$ be
vertices of $T$ adjacent to $u$ in $G$ , if $x_{1}^{-}=x_{2}^{-}=z$
for some vertex $z$ of $T$, let $T^{\prime }=\left( T+ux_{1}\right)
-zx_{1}$. Since $zx_{1}$ and $zx_{2}$ are edges of $ T$,
$d_{T^{\prime }}\left( z\right) $ $\geq 2$ and $T^{\prime }$ is a
subtree of $G$ with $V\left( T^{\prime }\right) =V\left( T\right) $,
with less than $2+\sum\limits_{j=1}^{k}\left( d_{i_{j}}-2\right) $
leaves and such that $d_{T^{\prime }}\left( w_{j}\right) \leq d_{j}$
for each $w_{j}\in V\left( T^{\prime }\right) $. Again, this is a
contradiction, therefore $ \left\vert B\right\vert =d\left( u\right)
$.

Since no vertex in $A\cup \left( B\setminus \left\{ u\right\}
\right) $ is a leave of $T$,

\begin{center}
$\left\vert A\cup B\right\vert \leq \left\vert V\left( T\right)
\right\vert -n_{1}+1\leq \left( n-1\right)
-n_{1}+1=n-2-\sum\limits_{j=1}^{k}\left( d_{j}-2\right) $.
\end{center}

Also

\begin{center}
$\left\vert A\cup B\right\vert =\left\vert A\right\vert +\left\vert
B\right\vert -\left\vert A\cap B\right\vert =d\left( u\right)
+d\left( v\right) -\left\vert A\cap B\right\vert \geq
n-1-\sum\limits_{j=1}^{k}\left( d_{j}-2\right) -\left\vert A\cap
B\right\vert $.
\end{center}

Therefore $\left\vert A\cap B\right\vert \geq 1$; let $z^{-}\in
A\cap B$. We consider two cases:

\textit{Case} 1. Edge $z^{-}z$ lies on the path $T_{uv}$.

\begin{figure}[h!!]
\begin{center}
\includegraphics[height=3cm, width=7cm]{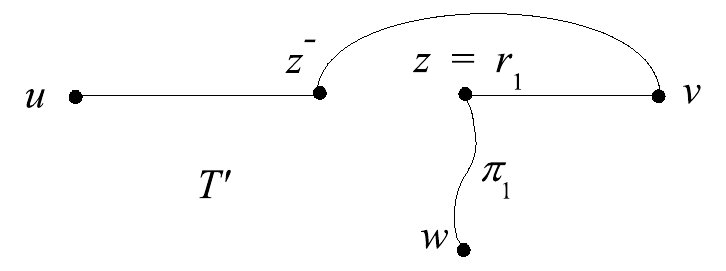}
\caption{$T^{\prime }=\left( \left( T+z^{-}v\right) -z^{-}z\right)
\cup \pi _{1}$}
\label{two}
\end{center}
\end{figure}

If $z=r_{1}$(see Fig. \ref{two}), let

\[
T^{\prime }=\left( \left( T+z^{-}v\right) -z^{-}z\right) \cup \pi
_{1}\text{ and}
\]

\begin{figure}[h!!]
\begin{center}
\includegraphics[height=3cm, width=8cm]{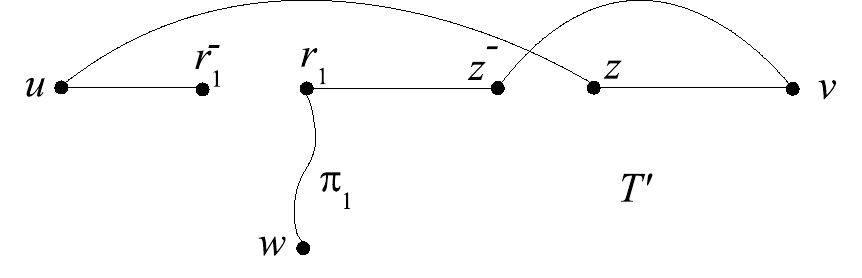}
\caption{$T^{\prime }=\left( \left( \left( \left( T+uz\right)
+z^{-}v\right) -r_{1}^{-}r_{1}\right) -z^{-}z\right) \cup \pi _{1}$}
\label{three}
\end{center}
\end{figure}

and if $r_{1}\neq z$ (see Fig. \ref{three}), let
\[
T^{\prime }=\left( \left( \left( \left( T+uz\right) +z^{-}v\right)
-r_{1}^{-}r_{1}\right) -z^{-}z\right) \cup \pi _{1}\text{.}
\]

Both situations lead to a contradiction since $T^{\prime }$ is a
subtree of $ G$ larger than $T$, with at most
$2+\sum\limits_{j=1}^{k}\left( d_{j}-2\right) $ leaves and such that
$d_{T}\left( w_{j}\right) \leq d_{j}$ for each $w_{j}\in V\left(
T^{\prime }\right) $.

\textit{Case} 2. Edge $z^{-}z$ does not lie on the path $T_{uv}$.

\begin{figure}[h!!]
\begin{center}
\includegraphics[height=3cm, width=8cm]{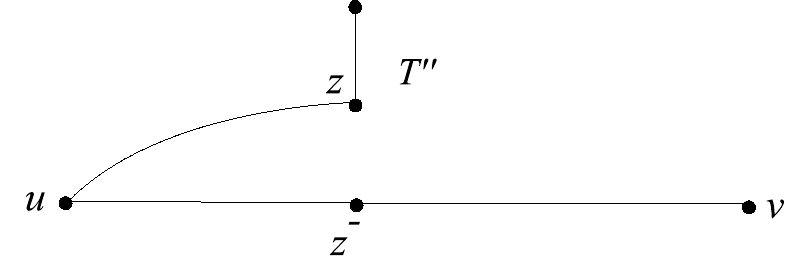}
\caption{$T^{\prime \prime }=\left( T+uz\right) -z^{-}z$}
\label{four}
\end{center}
\end{figure}

If $z^{-}$lies in $T_{uv}$, let $T^{\prime \prime }=\left(
T+uz\right) -z^{-}z$ (see Fig. \ref{four}).

\begin{figure}[h!!]
\begin{center}
\includegraphics[height=4cm, width=8cm]{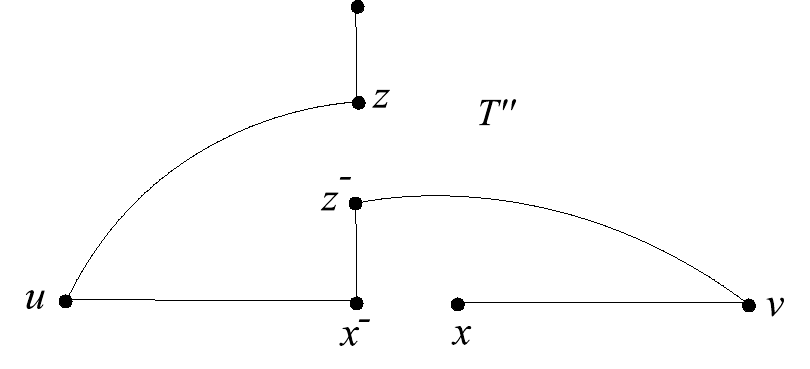}
\caption{$T^{\prime \prime }=\left( \left( \left( \left( T+uz\right)
+z^{-}v\right) -x^{-}x\right) -z^{-}z\right)$}
\label{five}
\end{center}
\end{figure}

And if $z^{-}$ does not lie in $T_{uv}$, let $x$ be a vertex in
$T_{uv}$ not in $T_{uz^{-}}$ such that $x^{-}$ is a vertex in
$T_{uz^{-}}$ (see Fig. \ref {five}). Let
\[
T^{\prime \prime }=\left( \left( \left( \left( T+uz\right)
+z^{-}v\right) -x^{-}x\right) -z^{-}z\right)
\]

In this case $T^{\prime \prime }$ is a subtree of $G$ with $V\left(
T^{\prime \prime }\right) =V\left( T\right) $, with at most
$n_{1}-1=$ $ 1+\sum\limits_{j=1}^{k}\left( d_{j}-2\right) $ leaves
and such that $ d_{T^{\prime \prime }}\left( w_{j}\right) \leq
d_{j}$ for each $w_{j}\in V\left( T^{\prime \prime }\right) $. As
seen above, this is not possible.

Cases 1 and 2 cover all possibilities, therefore $T$ is a spanning
tree of $G$.

\bigskip

Let $k\geqslant 1$ and $d_{1},d_{2},\ldots ,d_{n}$ be integers with
$3\leq d_{1}\leq d_{2}\leq \cdots \leq d_{n}$ and $X=\left\{
x_{1},x_{2},\ldots ,x_{k}\right\} $ and $Y=\left\{
y_{1},y_{2},\ldots ,y_{2-k+d_{1}+\cdots +d_{k}}\right\} $ be sets of
vertices. The complete bipartite graph$\ G$
with bipartition $\left( X,Y\right) $ is $k$-connected, has $
n=2+\sum\limits_{j=1}^{k}d_{i}$ vertices and is such that $d\left(
u\right) +d\left( v\right) \geq 2k=n-2-\sum\limits_{j=1}^{k}\left(
d_{i}-2\right) $ for any vertices $u$ and $v$ of $G$. Nevertheless,
if $T$ is a spanning tree of $G$, then $d_{T}\left( x_{j}\right)
>d_{j}$ for some $j=1,2,\ldots ,k$. This shows that the condition in
Theorem 2 is tight.


\end{document}